\documentclass[12pt, reqno]{amsart}
\usepackage{amsmath, amsthm, amscd, amsfonts, amssymb, graphicx, color}
\usepackage[bookmarksnumbered, colorlinks, plainpages]{hyperref}

\newtheorem{theorem}{Theorem}[section]
\newtheorem{lemma}[theorem]{Lemma}

\theoremstyle{definition}

\theoremstyle{remark}
\newtheorem{remark}[theorem]{Remark}
\numberwithin{equation}{section}

\begin{document}
\setcounter{page}{1}

\title[2-Local derivations on matrix rings]{2-Local derivations on matrix rings over associative
rings}

\author[Sh.\ Ayupov]{Shavkat Ayupov$^1$}

\address{$^{1}$ Department of Mathematics,
Institute of Mathematics, National University of Uzbekistan,
Tashkent, Uzbekistan.}
\email{\textcolor[rgb]{0.00,0.00,0.84}{sh$_-$ayupov@mail.ru}}

\author[F. Arzikulov]{Farhodjon Arzikulov$^2$}

\address{$^{2}$ Department of Mathematics, Andizhan State University, Andizhan 710020, Uzbekistan.}
\email{\textcolor[rgb]{0.00,0.00,0.84}{arzikulovfn@rambler.ru;
arzikulovf@mail.ru}}




\thanks{Supported by TWAS, The Abdus Salam, International Centre, for
Theoretical Physics (ICTP), Grant: 13-244RG/MATHS/AS$_-$I-UNESCO
FR:3240277696}

\begin{abstract}
In the present paper it is proved that every inner 2-local
derivation on the matrix ring $M_n(\Re)$ of $n\times n$ matrices over a commutative
associative ring $\Re$ is an inner derivation. Also, it is proved that,
every derivation on an associative ring $\Re$ has an extension to a derivation on
the matrix ring $M_n(\Re)$ of $n\times n$ matrices over $\Re$.
\end{abstract}

\maketitle

\section{Introduction}

The present paper is devoted to 2-local derivations on associative
rings. Recall that a 2-local derivation is defined as follows:
given a ring $\Re$, a map $\Delta : \Re \to \Re$ (not additive in
general) is called a 2-local derivation if for every $x$, $y\in
\Re$, there exists a derivation $D_{x,y} : \Re\to \Re$ such that
$\Delta(x)=D_{x,y}(x)$ and $\Delta(y)=D_{x,y}(y)$.

In 1997, P. \v{S}emrl \cite{S} introduced the notion of 2-local
derivations and described 2-local derivations on the algebra
$B(H)$ of all bounded linear operators on the infinite-dimensional
separable Hilbert space H. A similar description for the
finite-dimensional case appeared later in \cite{KK}. In the paper
\cite{LW} 2-local derivations have been described on matrix
algebras over finite-dimensional division rings. In \cite{AK} the
authors suggested a new technique and have generalized the above
mentioned results of \cite{S} and \cite{KK} for arbitrary Hilbert
spaces. Namely they considered 2-local derivations on the algebra
$B(H)$ of all linear bounded operators on an arbitrary (no
separability is assumed) Hilbert space $H$ and proved that every
2-local derivation on $B(H)$ is a derivation. In \cite{AA2},
\cite{AK2} the authors extended the above results and give a proof
of the theorem for arbitrary von Neumann algebras.

In this article we develop an algebraic approach to the
investigation of derivations and 2-local derivations on
associative rings. Since we consider a sufficiently general case
of associative rings we restrict our attention only on inner
derivations and inner 2-local derivations. In particular, we
consider the following problem:  if an inner 2-local derivation on
an associative ring is a derivation then is the latter derivation
inner? The answer to this question is affirmative if the ring is
generated by two elements (Theorem \ref{10}).

In this article  we consider 2-local derivations on the matrix
ring $M_n(\Re)$ over an associative ring $\Re$. The first step of
the investigation consists of proving that, in the case of a
commutative associative ring $\Re$ an arbitrary inner 2-local
derivation on $M_n(\Re)$ is an inner derivation. This result
extends the one obtained in \cite{LW} to the infinite dimensional
case but for a commutative ring $\Re$.

The second step consists of showing that every derivation on an
associative ring $\Re$ has an extension to a derivation on
the matrix ring $M_n(\Re)$ of $n\times n$ matrices over $\Re$.

\section{2-local derivations on matrix rings}

\medskip

Let $\Re$ be a ring. Recall that a map $D : \Re\to \Re$ is called
a derivation, if $D(x+y)=D(x)+D(y)$ and $D(xy)=D(x)y+xD(y)$ for
any two elements $x$, $y\in \Re$. A derivation $D$ on a ring $\Re$
is called an inner derivation, if there exists an element $a\in
\Re$ such that
$$
D(x)=ax-xa, x\in \Re.
$$

A map $\Delta : \Re\to \Re$ is called a 2-local derivation, if for
any two elements $x$, $y\in \Re$ there exists a derivation
$D_{x,y}:\Re\to \Re$ such that $\Delta (x)=D_{x,y}(x)$, $\Delta
(y)=D_{x,y}(y)$.

A map $\Delta : \Re\to \Re$ is called an inner 2-local derivation,
if for any two elements $x$, $y\in \Re$ there exists an element
$a\in \Re$ such that $\Delta (x)=ax-xa$, $\Delta (y)=ay-ya$.

Let $\Re$ be an associative ring with identity, $M_n(\Re)$ be the
matrix ring over $\Re$, $n>1$, of matrices of the form
$$
\left[%
\begin{array}{cccc}
 a^{1,1} & a^{1,2} & \cdots & a^{1,n}\\
a^{2,1} & a^{2,2} & \cdots & a^{2,n}\\
\vdots & \vdots & \ddots & \vdots\\
a^{n,1} & a^{n,2} & \cdots & a^{n,n}\\
\end{array}%
\right],
a^{i,j}\in \Re, i,j=1,2,\dots ,n.
$$
Let $\{e_{i,j}\}_{i,j=1}^n$ be the set of matrix units in
$M_n(\Re)$, i.e. $e_{i,j}$ is a matrix with components
$a^{i,j}={\bf 1}$ and $a^{k,l}={\bf 0}$ if $(i,j)\neq(k,l)$, where
${\bf 1}$ is an identity element, ${\bf 0}$ is the zero element of
$\Re$ and a matrix $a\in M_n(\Re)$ is written as $a=\sum_{k,l=1}^n
a^{k,l}e_{k,l}$, where $a^{k,l}\in \Re$ for $k,l=1,2,\dots, n$.
Let $\bar{M}_2(\Re)$ be the subring of $M_n(\Re)$, generated by
the subsets $\Re e_{i,j}$, $i,j=1,2$ in $M_n(\Re)$. It is clear
that
$$
\bar{M}_2(\Re)\cong M_2(\Re).
$$

The following theorem is the main result of the paper.

\begin{theorem} \label{1}
Let $\Re$ be a commutative associative ring with identity, and let $M_n(\Re)$
be the matrix ring of $n\times n$ matrices over $\Re$, $n>1$. Then
any inner 2-local derivation on the matrix ring $M_n(\Re)$ is an inner derivation.
\end{theorem}

\medskip

First let us prove lemmata and a theorem which will be used in the
proof of theorem \ref{1}. Throughout, $\Re$ denotes an associative
ring with identity, $M_n(\Re)$ denotes the matrix ring of $n\times
n$ matrices over $\Re$, $n>1$. Let $\Delta :M_n(\Re)\to M_n(\Re)$
be an inner 2-local derivation. Consider the subset
$\{a(i,j)\}_{i,j=1}^n\subset M_n(\Re)$  such that
$$
\Delta(e_{i,j})=a(i,j)e_{i,j}-e_{i,j}a(i,j),
$$
$$
\Delta(\sum_{k=1}^{n-1}e_{k,k+1})=a(i,j)(\sum_{k=1}^{n-1}e_{k,k+1})-(\sum_{k=1}^{n-1}e_{k,k+1})a(i,j).
$$
Put $a_{i,j}=e_{i,i}a(j,i)e_{j,j}$, for all pairs of different
indices $i$, $j$ and let $\sum_{k\neq l} a_{k,l}$ be the sum of
all such elements.

\begin{lemma} \label{21}
Let $\Delta :M_n(\Re)\to M_n(\Re)$ be an inner 2-local derivation.
Then the following equalities hold for any pair $i$, $k$ of
different indices in $\{1,2,\dots,n\}$
$$
e_{k,k}a(i,j)e_{i,j}=e_{k,k}a(i,k)e_{i,j}
$$
and for any pair $j$, $k$ of different indices in
$\{1,2,\dots,n\}$
$$
e_{i,j}a(i,j)e_{k,k}=e_{i,j}a(k,j)e_{k,k}.
$$
\end{lemma}

\begin{proof}
Let $d\in M_n(\Re)$ be such element that
$$
\Delta(e_{i,j})=de_{i,j}-e_{i,j}d,
\Delta(e_{i,k})=de_{i,k}-e_{i,k}d.
$$
Then
$$
a(i,j)e_{i,j}-e_{i,j}a(i,j)=de_{i,j}-e_{i,j}d,
$$
$$
a(i,k)e_{i,k}-e_{i,k}a(i,k)=de_{i,k}-e_{i,k}d,
$$
and
$$
e_{k,k}a(i,j)e_{i,j}=e_{k,k}de_{i,j},
e_{k,k}a(i,k)e_{i,k}=e_{k,k}de_{i,k}.
$$
Hence
$$
e_{k,k}a(i,k)e_{i,k}e_{k,j}=e_{k,k}de_{i,k}e_{k,j}=e_{k,k}de_{i,j}
$$
and
$$
e_{k,k}a(i,k)e_{i,j}=e_{k,k}de_{i,j}.
$$
Similarly,
$$
e_{i,j}a(i,j)e_{k,k}=e_{i,j}a(k,j)e_{k,k}.
$$
\end{proof}

\medskip

\begin{lemma} \label{2}
Let $\Delta :M_n(\Re)\to M_n(\Re)$ be an inner 2-local derivation.
Then for any pair $i$, $j$ of different indices in
$\{1,2,\dots,n\}$ the following equality holds
$$
\Delta(e_{i,j})=(\sum_{k,l=1, k\neq
l}^na_{k,l})e_{i,j}-e_{i,j}(\sum_{k,l=1, k\neq
l}^na_{l,k})+a(i,j)^{i,i}e_{i,j}-e_{i,j}a(i,j)^{j,j},
\,\,\,\,\,\,\,(1)
$$
where $a(i,j)_{i,i}$, $a(i,j)_{j,j}$ are the appropriate
components of the matrices $e_{i,i}a(i,j)e_{i,i}$,
$e_{j,j}a(i,j)e_{j,j}$.
\end{lemma}

\begin{proof}
We have
$$
\Delta(e_{i,j})=a(i,j)e_{i,j}-e_{i,j}a(i,j)=\sum_{k=1}^ne_{k,k}
a(i,j)e_{i,j}-\sum_{k=1}^ne_{i,j}a(i,j)e_{k,k}=
$$
$$
\sum_{k=1, k\neq i}^ne_{k,k} a(i,j)e_{i,j}-\sum_{k=1, k\neq
j}^ne_{i,j}a(i,j)e_{k,k}+e_{i,i}
a(i,j)e_{i,j}+e_{i,j}a(i,j)e_{j,j}=
$$
$$
\sum_{k=1, k\neq i}^ne_{k,k} a(i,k)e_{i,j}-\sum_{k=1, k\neq
j}^ne_{i,j}a(k,j)e_{k,k}+a(i,j)^{i,i}e_{i,j}-e_{i,j}a(i,j)^{j,j}=
$$
$$
\sum_{k=1, k\neq i}^na_{k,i}e_{i,j}-\sum_{k=1, k\neq
j}^ne_{i,j}a_{j,k}+a(i,j)^{i,i}e_{i,j}-e_{i,j}a(i,j)^{j,j}=
$$
$$
(\sum_{k,l=1, k\neq l}^na_{k,l})e_{i,j}-e_{i,j}(\sum_{k,l=1, k\neq
l}^na_{l,k})+a(i,j)^{i,i}e_{i,j}-e_{i,j}a(i,j)^{j,j}
$$
by lemma \ref{21}.
\end{proof}

Let $x_o=\sum_{k=1}^{n-1}e_{k,k+1}$. Then there exists an element
$c\in M_n(\Re)$ such that
$$
\Delta(x_o)=cx_o-x_oc.
$$
Let $c=\sum_{i,j=1}^nc_{i,j}$ be the Pierce decomposition of $c$,
where $c_{i,j}=e_{i,i}ce_{j,j}$, $i,j=1,2,\dots, n$.

\begin{lemma} \label{3}
Let $\Delta :M_n(\Re)\to M_n(\Re)$ be an inner 2-local derivation.
Let $k$, $l$ be arbitrary different numbers in $\{1,2,\dots,n\}$,
$b\in M_n(\Re)$ be an element such that
$$
\Delta(x_o)=bx_o-x_ob.
$$
Then $c^{k,k}-c^{l,l}=b^{k,k}-b^{l,l}$, where
$c_{i,i}=c^{i,i}e_{i,i}$, $b_{i,i}=b^{i,i}e_{i,i}$, $c^{i,i}$,
$b^{i,i}\in \Re$, $i=1,2,\dots,n$.
\end{lemma}

\begin{proof} We may assume that  $k<l$. We have
$$
\Delta(x_o)=cx_o-x_oc=bx_o-x_ob.
$$
Hence
$$
e_{k,k}(cx_o-x_oc)e_{k+1,k+1}=e_{k,k}(bx_o-x_ob)e_{k+1,k+1}
$$
and
$$
c^{k,k}-c^{k+1,k+1}=b^{k,k}-b^{k+1,k+1}.
$$
Then for the sequence
$$
(k,k+1),(k+1,k+2)\dots (l-1,l)
$$
we have
$$
c^{k,k}-c^{k+1,k+1}=b^{k,k}-b^{k+1,k+1},
c^{k+1,k+1}-c^{k+2,k+2}=b^{k+1,k+1}-b^{k+2,k+2},\dots
$$
$$
c^{l-1,l-1}-c^{l,l}=b^{l-1,l-1}-b^{l,l}.
$$
Hence
$$
c^{k,k}-b^{k,k}=c^{k+1,k+1}-b^{k+1,k+1},
c^{k+1,k+1}-b^{k+1,k+1}=c^{k+2,k+2}-b^{k+2,k+2},\dots
$$
$$
c^{l-1,l-1}-b^{l-1,l-1}=c^{l,l}-b^{l,l}.
$$
Therefore $c^{k,k}-b^{k,k}=c^{l,l}-b^{l,l}$, i.e.
$c^{k,k}-c^{l,l}=b^{k,k}-b^{l,l}$. The proof is complete.
\end{proof}

\medskip

Let $a_{i,i}=c_{i,i}$ for $i=,2,\dots,n$ and
$\bar{a}=\sum_{i,j=1}^na_{i,j}$.

{\em A proof of theorem \ref{1}.}  
Let $\Delta :M_n(\Re)\to M_n(\Re)$ be an inner 2-local derivation,
$x$ be an arbitrary matrix in $M_n(\Re)$ and let $d(i,j)\in
M_n(\Re)$ be an element such that
$$
\Delta(e_{i,j})=d(i,j)e_{i,j}-e_{i,j}d(i,j), \,\,\,\,
\Delta(x)=d(i,j)x-xd(i,j)
$$
and $i\neq j$. Then by Lemma \ref{2} we have
$$
\Delta(e_{i,j})=d(i,j)e_{i,j}-e_{i,j}d(i,j)=
$$
$$
e_{i,i}d(i,j)e_{i,j}-e_{i,j}d(i,j)e_{j,j}+
(1-e_{i,i})d(i,j)e_{i,j}-e_{i,j}d(i,j)(1-e_{j,j})=
$$
$$
a(i,j)_{i,i}e_{i,j}-e_{i,j}a(i,j)_{j,j}+(\sum_{k\neq l}
a_{k,l})e_{i,j}-e_{i,j}(\sum_{k\neq l} a_{k,l})
$$
for all $i$, $j$ in $\{1,2,\dots,n\}$.

Since
$e_{i,i}d(i,j)e_{i,j}-e_{i,j}d(i,j)e_{j,j}=a(i,j)_{i,i}e_{i,j}-e_{i,j}a(i,j)_{j,j}$
we have
$$
(1-e_{i,i})d(i,j)e_{i,i}=(\sum_{k\neq l} a_{k,l})e_{i,i},
$$
$$
e_{j,j}d(i,j)(1-e_{j,j})=e_{j,j}(\sum_{k\neq l} a_{k,l})
$$
for all different numbers $i$ and $j$ in $\{1,2,\dots,n\}$.

Hence
$$
e_{j,j}\Delta(x)e_{i,i}=e_{j,j}(d(i,j)x-xd(i,j))e_{i,i}=
$$
$$
e_{j,j}d(i,j)(1-e_{j,j})xe_{i,i}+
e_{j,j}d(i,j)e_{j,j}xe_{i,i}-e_{j,j}x(1-e_{i,i})d(i,j)e_{i,i}-e_{j,j}xe_{i,i}d(i,j)e_{i,i}=
$$
$$
e_{j,j}(\sum_{k\neq l} a_{k,l})xe_{i,i}-e_{j,j}x(\sum_{k\neq l}
a_{k,l})e_{i,i}+
e_{j,j}d(i,j)e_{j,j}xe_{i,i}-e_{j,j}xe_{i,i}d(i,j)e_{i,i}.
$$
We have
$$
\Delta(\sum_{k=1}^{n-1}e_{k,k+1})=a(i,j)(\sum_{k=1}^{n-1}e_{k,k+1})-(\sum_{k=1}^{n-1}e_{k,k+1})a(i,j)
$$
by the definition of $a(i,j)$. Then by lemma \ref{3} we have
$$
a(i,j)^{j,j}-a(i,j)^{i,i}=c^{j,j}-c^{i,i},
$$
where
$$
c_{k,k}=c^{k,k}e_{k,k}, c^{k,k}\in \Re, k=1,2,\dots,n,
$$
$$
a(i,j)=\sum_{kl=1}^n a(i,j)^{k,l}e_{k,l}, a(i,j)^{k,l}\in \Re,
k,l=1,2,\dots,n.
$$
Since
$$
d(i,j)e_{i,j}-e_{i,j}d(i,j)=a(i,j)e_{i,j}-e_{i,j}a(i,j)
$$
we have
$$
e_{i,i}d(i,j)e_{i,j}-e_{i,j}d(i,j)e_{j,j}=e_{i,i}a(i,j)e_{i,j}-e_{i,j}a(i,j)e_{j,j}
$$
and
$$
(d(i,j)^{i,i}-d(i,j)^{j,j})e_{i,j}=(a(i,j)^{i,i}-a(i,j)^{j,j})e_{i,j},
$$
where $d(i,j)=\sum_{kl=1}^n d(i,j)^{k,l}e_{k,l}$. Hence
$d(i,j)^{i,i}-d(i,j)^{j,j}=a(i,j)^{i,i}-a(i,j)^{j,j}$, i.e.
$d(i,j)^{j,j}-d(i,j)^{i,i}=a(i,j)^{j,j}-a(i,j)^{i,i}$. Therefore
$$
e_{j,j}d(i,j)e_{j,j}xe_{i,i}-e_{j,j}xe_{i,i}d(i,j)e_{i,i}=d(i,j)^{j,j}x^{j,i}e_{i,j}-x^{j,i}d(i,j)^{i,i}e_{j,i}=
$$
$$
(d(i,j)^{j,j}-d(i,j)^{i,i})x^{j,i}e_{j,i}=(a(i,j)^{j,j}-a(i,j)^{i,i})x^{j,i}e_{j,i}=
$$
$$
(c^{j,j}-c^{i,i})x^{j,i}e_{j,i}=c^{j,j}x^{j,i}e_{j,i}-x^{j,i}c^{i,i}e_{j,i}=(c^{j,j}e_{j,j})e_{j,j}xe_{i,i}-e_{j,j}xe_{i,i}(c^{i,i}e_{i,i})=
$$
$$
a_{j,j}e_{j,j}xe_{i,i}-e_{j,j}xe_{i,i}a_{i,i},
$$
where $x=\sum_{kl=1}^n x^{k,l}e_{k,l}$. Hence
$$
e_{j,j}\Delta(x)e_{i,i}=e_{j,j}(\sum_{k\neq l}
a_{k,l})xe_{i,i}-e_{j,j}x(\sum_{k\neq l} a_{k,l})e_{i,i}+
a_{j,j}e_{j,j}xe_{i,i}-e_{j,j}xe_{i,i}a_{i,i}=
$$
$$
e_{j,j}(\sum_{k\neq l} a_{k,l})xe_{i,i}-e_{j,j}x(\sum_{k\neq l}
a_{k,l})e_{i,i}+e_{j,j}(\sum_{k=1}^na_{k,k})xe_{i,i}-e_{j,j}x(\sum_{k=1}^na_{k,k})e_{i,i}=
$$
$$
e_{j,j}(\sum_{kl=1}^n a_{k,l})xe_{i,i}-e_{j,j}x(\sum_{kl=1}^n
a_{k,l})e_{i,i}=e_{j,j}(\bar{a}x-x\bar{a})e_{i,i}.
$$

Let $d(i,i)$, $v$, $w\in \mathcal{M}$ be elements such that
$$
\bigtriangleup(e_{i,i})=d(i,i)e_{i,i}-e_{i,i}d(i,i), \,\,\,\,
\bigtriangleup(x)=d(i,i)x-xd(i,i),
$$
$$
\bigtriangleup(e_{i,i})=ve_{i,i}-e_{i,i}v,
\bigtriangleup(e_{i,j})=ve_{i,j}-e_{i,j}v,
$$
and
$$
\bigtriangleup(e_{i,i})=we_{i,i}-e_{i,i}w,
\bigtriangleup(e_{j,i})=we_{j,i}-e_{j,i}w.
$$
Then
$$
(1-e_{i,i})a(i,j)e_{i,i}=(1-e_{i,i})ve_{i,i}=(1-e_{i,i})d(i,i)e_{i,i},
$$
and
$$
e_{i,i}a(j,i)(1-e_{i,i})=e_{i,i}w(1-e_{i,i})=e_{i,i}d(i,i)(1-e_{i,i}).
$$
By lemma \ref{2} we have
$$
\Delta(e_{i,j})=a(i,j)e_{i,j}-e_{i,j}a(i,j)=
$$
$$
(\sum_{k\neq l} a_{k,l})e_{i,j}-e_{i,j}(\sum_{k\neq l}
a_{k,l})+a(i,j)_{i,i}e_{i,j}-e_{i,j}a(i,j)_{j,j}
$$
and
$$
(1-e_{i,i})a(i,j)e_{i,i}=(\sum_{k\neq l} a_{k,l})e_{i,i}.
$$
Similarly
$$
e_{i,i}a(j,i)(1-e_{i,i})=e_{i,i}(\sum_{k\neq l} a_{k,l}).
$$

Hence
$$
e_{i,i}\Delta(x)e_{i,i}=e_{i,i}(d(i,i)x-xd(i,i))e_{i,i}=
$$
$$
e_{i,i}d(i,i)(1-e_{i,i})xe_{i,i}+
e_{i,i}d(i,i)e_{i,i}xe_{i,i}-e_{i,i}x(1-e_{i,i})d(i,i)e_{i,i}-e_{i,i}xe_{i,i}d(i,i)e_{i,i}=
$$
$$
e_{i,i}a(j,i)(1-e_{i,i})xe_{i,i}+
e_{i,i}d(i,i)e_{i,i}xe_{i,i}-e_{i,i}x(1-e_{i,i})a(i,j)e_{i,i}-e_{i,i}xe_{i,i}d(i,i)e_{i,i}=
$$
$$
e_{i,i}(\sum_{k\neq l} a_{k,l})xe_{i,i}-e_{i,i}x(\sum_{k\neq l}
a_{k,l})e_{i,i}+
e_{i,i}d(i,i)e_{i,i}xe_{i,i}-e_{i,i}xe_{i,i}d(i,i)e_{i,i}=
$$
$$
e_{i,i}(\sum_{k\neq l} a_{k,l})xe_{i,i}-e_{i,i}x(\sum_{k\neq l}
a_{k,l})e_{i,i}+0=
$$
$$
e_{i,i}(\sum_{k\neq l} a_{k,l})xe_{i,i}-e_{i,i}x(\sum_{k\neq l}
a_{k,l})e_{i,i}+c_{i,i}e_{i,i}xe_{i,i}-e_{i,i}xc_{i,i}e_{i,i}=
$$
$$
e_{i,i}(\sum_{k\neq l} a_{k,l})xe_{i,i}-e_{i,i}x(\sum_{k\neq l}
a_{k,l})e_{i,i}+
$$
$$
e_{i,i}(\sum_{k=1}^na_{k,k})xe_{i,i}-e_{i,i}x(\sum_{k=1}^na_{k,k})e_{i,i}=
$$
$$
e_{i,i}(\sum_{kl=1}^n a_{k,l})xe_{i,i}-e_{i,i}x(\sum_{kl=1}^n
a_{k,l})e_{i,i}=e_{i,i}(\bar{a}x-x\bar{a})e_{i,i}.
$$

By the above conclusions we have
$$
\Delta(x)=\sum_{kl=1}^ne_{k,k}\Delta(x)e_{l,l}=\sum_{kl=1}^ne_{k,k}(\bar{a}x-x\bar{a})e_{l,l}=\bar{a}x-x\bar{a}
$$
for all $x\in M_n(\Re)$. The proof is complete.
$\triangleright$

\bigskip

\section{On extensions of derivations and 2-local
derivations}

\begin{lemma} \label{5}
Let $M_2(\Re)$ be the matrix ring of $2\times 2$ matrices over an
associative ring $\Re$ with identity and let $D$ be a derivation
on the subring $\Re e_{1,1}$ and $\delta$ be a derivation on $\Re$
induced by $D$. Then the map
$$
\bar{D}
\left(\left[%
\begin{array}{lr}
  \lambda & \mu \\
  \nu & \eta \\
\end{array}%
\right]\right)=
\left[%
\begin{array}{lr}
  \delta(\lambda) & \mu \\
  -\nu & \delta(\eta) \\
\end{array}%
\right], \lambda, \mu, \nu, \eta\in \Re,
$$
is a derivation.
\end{lemma}

{\em Proof} It is easy to check that for $a,b\in M_2(\Re)$ we have
$\bar{D}(ab)=\bar{D}(a)b+a\bar{D}(b)$. Indeed, the map $\bar{D}$
is just equal to $\bar{\delta}+d_U$, where
$$
\bar{\delta}\left(\left[%
\begin{array}{lr}
\lambda & \mu \\
\nu & \eta \\
\end{array}%
\right]\right)=
\left[%
\begin{array}{lr}
\delta(\lambda) & 0 \\
  0 & \delta(\eta) \\
\end{array}%
\right], \lambda, \mu, \nu, \eta\in \Re,
$$
and $d_U$ is the inner derivation induced by the matrix
$$
\left[
\begin{array}{lr} \frac{1}{2}  & 0 \\ 0 & -\frac{1}{2}
\end{array} \right].
$$
$\triangleright$

Let $\bar{M}_m(\Re)$ be a subring of $M_n(\Re)$, $m<n$, generated
by the subsets
$$
\Re e_{i,j}, i,j=1,2,\dots, m
$$
in $M_n(\Re)$. It is clear that
$$
\bar{M}_m(\Re)\cong M_m(\Re).
$$

\begin{lemma} \label{6}
Let $\Re$ be an associative ring, and let $M_n(\Re)$ be a matrix
ring of $n\times n$ matrices over $\Re$, $n>2$. Then every
derivation on $\bar{M}_2(\Re)$ can be extended to a derivation on
$M_n(\Re)$.
\end{lemma}

\begin{proof} By lemma \ref{5} every derivation on $\bar{M}_2(\Re)$
can be extended to a derivation on $M_4(\Re)$. In its turn, every
derivation on $\bar{M}_4(\Re)$ can be extended to a derivation on
$M_8(\Re)$ and so on. Thus every derivation $\partial$ on
$\bar{M}_2(\Re)$ can be extended to a derivation $D$ on
$M_{2^k}(\Re)$. Suppose that $n\leq 2^k$. Let $e=\sum_{i=1}^n
e_{i,i}$ and
$$
\bar{D}(a)=eD(a)e, a\in \bar{M}_n(\Re).
$$
Then $\bar{D}:\bar{M}_n(\Re)\to \bar{M}_n(\Re)$ and $\bar{D}$ is a
derivation on $\bar{M}_n(\Re)$. Indeed, it is clear that $\bar{D}$
is a linear map. At the same time, for all $a$, $b\in
\bar{M}_n(\Re)$ we have
$$
\bar{D}(ab)=eD(ab)e=e(D(a)b+aD(b))e=
$$
$$
eD(a)be+eaD(b)e=eD(a)eb+aeD(b)e=\bar{D}(a)b+a\bar{D}(b).
$$
Hence, $\bar{D}$ is a derivation. At the same time, the derivation
$\bar{D}$ coincides with the derivation $\partial$ on
$\bar{M}_2(\Re)$. Therefore, $\bar{D}$ is an extension of
$\partial$ to $\bar{M}_n(\Re)$. Hence every derivation $\partial$
on $\bar{M}_2(\Re)$ can be extended to a derivation on $M_n(\Re)$.
\end{proof}

Thus, in the case of the ring $M_2(\Re)$ for any derivation on the
subring $\Re e_{1,1}$ we can take its extension onto the whole
$M_2(\Re)$ defined as in lemma \ref{5}, which is also a
derivation.

\begin{theorem} \label{7}
Let $\Re$ be an associative ring, and let $M_n(\Re)$ be a matrix
ring of $n\times n$ matrices over $\Re$, $n>2$. Then every
derivation on $\Re$ can be extended to a derivation on
$M_n(\Re)$.
\end{theorem}

\begin{proof}
Let $\delta$ be an arbitrary derivation on $\Re$ and $D$ be the derivation on the subring $\Re e_{1,1}$
such that $\delta$ is induced by $D$. By lemma \ref{5} every derivation on $\Re e_{1,1}$ has an extension
to a derivation on the matrix ring $\bar{M}_2(\Re)$ and every derivation on $\bar{M}_2(\Re)$ 
has an extension to a derivation on the matrix ring $M_n(\Re)$ by lemma \ref{6}. Thus the statement of the
theorem is valid.
\end{proof}

\medskip

\begin{remark}
If $\Re$ is an arbitrary associative ring with identity and $n=1$ then 
theorem \ref{1} is not valid. Indeed, in \cite{ZL} there is an
example of an inner 2-local derivation on the associative ring
$U_2({\Bbb C})$ of $2\times 2$ triangle matrices
$$
\left[%
\begin{array}{cc}
  a^{1,1} & a^{1,2} \\
  0 & a^{2,2} \\
\end{array}%
\right], a^{1,1}, a^{1,2}, a^{2,2}\in {\Bbb C}
$$
over ${\Bbb C}$ (the complex numbers), which is not additive.

Also by \cite[Theorem 3.5]{AK4}
the lattice $P(\mathcal{M})$ of projections in a von Neumann algebra
$\mathcal{M}$ is not atomic if and only if the algebra $S(\mathcal{M})$ of
all measurable operators affiliated with $\mathcal{M}$ admits a 2-local derivation
which is not a derivation. Hence, if $\Re$ is the algebra $S(\mathcal{M})$ and 
$P(\mathcal{M})$ is not atomic then by \cite[Theorem 4.3]{AK4} no any
2-local derivation on $\Re e_{1,1}$ has an extension to a 2-local derivation on
$M_n(\Re)$, $n>1$.
\end{remark}

\medskip

We conclude the paper by the following more general observation.

\begin{theorem} \label{10}
Let $\Delta :\Re\to \Re$ be an inner 2-local derivation on an
associative ring $\Re$. Suppose that $\Re$ is generated by its two
elements. Then, if $\Delta$ is additive then it is an inner
derivation.
\end{theorem}

\begin{proof}
Let $x$, $y$ be generators of $\Re$, i.e.
$\Re=Alg(\{x,y\})$, where $Alg(\{x,y\})$ is an associative ring,
generated by the elements $x$, $y$ in $\Re$. We have that there
exists $d\in \Re$ such that
$$
\Delta(x)=[d,x], \Delta(y)=[d,y],
$$
where $[d,a]=da-ad$ for any $a\in \Re$.

Hence by the additivity of $\Delta$ we have
$$
\Delta(x+y)=\Delta(x)+\Delta(y)=[d,x+y].
$$
Note that
$$
\Delta(xy)=\Delta(x)y+x\Delta(y)=[d,x]y+x[d,y]=[d,xy],
$$
$$
\Delta(x^2)=\Delta(x)x+x\Delta(x)=[d,x]x+x[d,x]=[d,x^2],
$$
$$
\Delta(y^2)=\Delta(y)y+y\Delta(y)=[d,y]y+y[d,y]=[d,y^2],
$$
Similarly
$$
\Delta(x^k)=[d,x^k], \Delta(y^m)=[d,y^m],
\Delta(x^ky^m)=[d,x^ky^m]
$$
and
$$
\Delta(x^ky^mx^l)=\Delta(x^ky^m)x^l+x^ky^m\Delta(x^l)=[d,x^ky^m]x^l+x^ky^m[d,x^l]=[d,x^ky^mx^l].
$$
Finally, for every polynomial $p(x_1,x_2,\dots,x_m)\in \Re$, where
$x_1,x_2,\dots,x_m\in \{x,y\}$ we have
$$
\Delta(p(x_1,x_2,\dots,x_m))=[d,p(x_1,x_2,\dots,x_m)],
$$
i.e. $\Delta$ is an inner derivation on $\Re$. \end{proof}


\end{document}